\documentclass{jhrs}

\usepackage{amsmath,amssymb,amscd,amsthm,latexsym}
\usepackage{amsfonts}
\usepackage[english]{babel}
\usepackage[latin1]{inputenc}
\usepackage{graphicx}
\usepackage{float}
\usepackage{subfigure}
\usepackage[all]{xy}

\newtheorem{theor}{Theorem}[section]

\newtheorem{prop}[theor]{Proposition}
\newtheorem{cor}[theor]{Corollary}
\newtheorem{defin}[theor]{Definition}
\newtheorem{rem}[theor]{Remark}

\newtheorem{ex}[theor]{Example}

\newcommand{\dem}{\noindent \emph{Proof. }}
\newcommand{\findem}{\hfill $\Box$\\}

\newdir{ >}{{}*!/-7pt/\dir{>}}

\begin{document}

\title{Some collapsing operations for 2-dimensional precubical sets}

\author{Thomas Kahl}             
\email{kahl@math.uminho.pt}       
\address{Universidade do Minho\\ 
         Departamento de Matem\'atica e Aplicações\\
         Campus de Gualtar\\
         4710-057 Braga\\
         Portugal}
\thanks{This research has been supported by FEDER funds through ``Programa Operacional Factores de Competitividade - COMPETE'' and by FCT -  \emph{Fundação para a Ciência e a Tecnologia} through projects Est-C/MAT/UI0013/2011 and PTDC/MAT/0938317/2008.}



\begin{abstract}
In this paper, we consider 2-dimensional precubical sets, which can be used to model systems of two concurrently executing processes. From the point of view of concurrency theory, two precubical sets can be considered equivalent if their geometric realizations have the same directed homotopy type relative to the extremal elements in the sense of P. Bubenik. We give easily verifiable conditions under which it is possible to reduce a 2-dimensional precubical set to an equivalent smaller one by collapsing an edge or eliminat\nolinebreak ing a square and one or two free faces. We also look at some simple standard examples in order to illustrate how our results can be used to construct small models of 2-dimensional precubical sets.
\end{abstract}

\received{Month Day, Year}   
\revised{Month Day, Year}    
\published{Month Day, Year}  
\submitted{}  

\volumeyear{} 
\volumenumber{}  
\issuenumber{}   

\startpage{1}     

\maketitle


\maketitle 


\section{Introduction}

\begin{figure}
\begin{center}
\subfigure[]
{ 
\includegraphics[height=2.5cm]{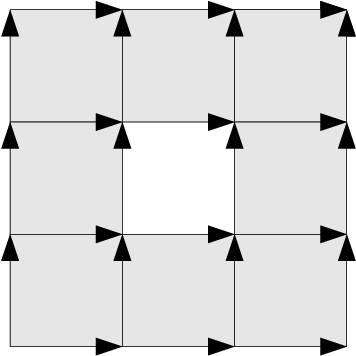}
\label{fig1}
}
\hspace{2cm}
\subfigure[]
{
\includegraphics[height=2.5cm]{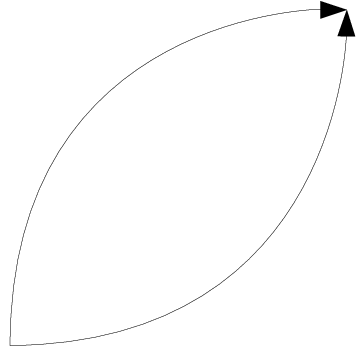}
\label{fig1b}
}
\caption{A very simple concurrent system and its fundamental bipartite graph}
\end{center}
\end{figure}

It has been known for some time now that precubical sets, i.e., cubical sets without degeneracies, can be used to model concurrent systems  (cf. \cite {Fajstrup}, \cite{FajstrupGR}, \cite{GaucherGoubault},  \cite{vanGlabbeek}, \cite{Goubault}, \cite{GoubaultTransCub}). These are systems of two or more computational processes which may communicate, share resources, and execute in parallel. Let us consider, as an example, a very simple concurrent system where two processes $A$ and $B$ write to a piece of shared memory. Each process performes a sequence of three actions: it accesses the memory, writes its name, and terminates. The processes execute simultaneously but cannot write to the memory at the same time. This situation can be modeled by the 2-dimensional precubical set  depicted in figure \ref{fig1} (the definition of precubical sets is recalled in \ref{precube}). The vertices represent the states of the system, the horizontal arrows represent the actions of process $A$, and the vertical arrows represent the actions of process $B$. Moreover, if it does not matter in which order an action of $A$ and an action of $B$ are executed and they may actually be performed concurrently, then this is indicated by a square linking the two pairs of arrows corresponding to a consecutive execution of the actions. The precubical set has a hole reflecting the fact that only one process can write its name to the memory at a time.
\newline \indent Any precubical set can be realized geometrically as a topological space and indeed even as a d-space in the sense of M. Grandis \cite{GrandisBook}. A d-space is a topological space with a distinguished  set of paths, called d-paths, which equip the space with a direction of time. The d-paths in the geometric realization of a precubical set are obtained by pasting together increasing paths on cubes. In the interpretation of a precubical set as a model of a concurrent system, the passage to the geometric realization adds all possible intermediate states of the system to the model. The d-paths represent complete or partial executions of the system.\\	 

Consider again our example concurrent system. There exists an infinite number of d-paths leading from the initial state in the lower left corner to the final state in the upper right corner. Computer scientifically, two such d-paths can be considered equivalent if they represent executions which produce the same result, i.e., executions where the processes write to the memory in the same order. Geometrically, this happens precisely when the d-paths turn around the hole on the same side. This leads to the following notion of equivalence of d-paths: Two d-paths $\alpha$ and $\beta$ from a point $x$ in a d-space to a point $y$ are said to be dihomotopic relative to $\{0,1\}$ if there exists a homotopy $H$ from $\alpha$ to $\beta$ such that each path $H(-,t)$ is a d-path from $x$ to $y$. In the example, there are two relative dihomotopy classes of d-paths leading from the initial to the final state of the system, corresponding to the two possible orders in which the processes can write to the memory.
\newline \indent An important tool in the study of the directed structure of a d-space is its fundamental category (cf. \cite{Goubault}, \cite{GrandisBook}). This is the directed analogue of the fundamental groupoid of a topological space. The objects of the fundamental category of a d-space are its points and the morphisms are the relative dihomotopy classes of d-paths. The fundamental category of a d-space is, of course, a huge object, and a main line of research in the area of directed algebraic topology is the development of methods to extract the essential information of the fundamental category (cf. \cite{BubenikExtremal}, \cite{FajstrupGR}, \cite{Components}, \cite{Components2}, \cite[I.3]{GrandisBook}, \cite{RaussenInvariants}). A basic construction in this context is P. Bubenik's fundamental bipartite graph of a d-space, which is the full subcategory of the fundamental category generated by the so-called extremal elements (cf. \cite{BubenikExtremal}). In the geometric realization of a precubical set, the extremal elements are the points which correspond to the vertices in which no edge begins or no edge ends (cf. \ref{extrlP}). The fundamental bipartite graph of a d-space model of a concurrent system represents the essential execution schedules of the system. In the case of our example system, the fundamental bipartite graph is indicated in figure \ref{fig1b}.\\

In this paper, we are concerned with the problem of reducing a given precubical set to an equivalent smaller one. This approach complements the strategy to replace the fundamental category by a smaller object containing the relevant directed information. Ultimately, it would be very useful to have an efficient reduction algorithm for precubical sets. Collapses of cubes are certainly potential steps of such an algorithm. Unfortunately, the collapsing of cubes is a non-trivial matter in directed topology since eliminating a cube with a free face from a precubical set may change its directed structure. Perhaps surprisingly, it appears that collapsing operations preserve significantly more structure in 2-dimensional precubical sets than in higher dimensional ones. This paper is devoted to collapsing operations for 2-dimensional precubical sets. Collapsing operations for higher dimensional precubical sets will be discussed in a forthcoming paper. As the referee has pointed out, one may consider the 2-dimensional case particularly interesting because the full subcategory of the fundamental category generated by the vertices 
is in many important situations the same for a precubical set and its 2-skeleton (cf. \cite{Fajstrup}). 
\newline \indent The collapsing operations of this paper preserve the dihomotopy type relative to the extremal elements of the geometric realization. The notion of dihomotopy equivalence we use here is based on a straightforward extension of the notion of dihomotopy of d-paths to d-maps, i.e., morphisms of d-spaces (cf. \ref{defdihomotopy}). The reader should note  that other concepts of directed homotopy equivalence have been defined in the literature (cf. \cite{GaucherGoubault}, \cite{RaussenInvariants}). The concept of dihomotopy equivalence used in this paper is quite strong, and we show in \ref{FBG} that dihomotopy equivalences relative to the extremal elements between geometric realizations of precubical sets preserve the fundamental bipartite graph and the homotopy type of any trace space (cf. \cite{RaussenInvariants}, \cite{FahrenbergRaussen}, \cite{RaussenTraceCub}) between extremal elements. Using the collapsing operations of this paper, the precubical set represented in figure \ref{fig1} can be reduced to a precubical set that looks exactly like the fundamental bipartite graph in figure \ref{fig1b}.	\\

Our results may be seen as related to those of \cite{GaucherGoubault}. In that paper, the authors define S- and T-homotopy equivalences and give some of the collapsing operations we consider here as examples of such equivalences. Using these particular collapsing operations only, it is possible to reduce the precubical set of figure \ref{fig1} to one that looks like the fundamental bipartite graph in figure \ref{fig1b} with one additional arrow coming into the initial vertex and one additional arrow going out of the final vertex.\\
 
The main results are contained in sections 5 and 6. In section 2, we collect some basic facts on precubical sets and d-spaces. Section 3 is devoted to dihomotopy. In section 4, we present a method to construct d-maps and dihomotopies on geometric realizations of precubical sets. The last two sections are devoted to examples and some final remarks.

\section{Precubical sets and $\rm d$-spaces} 

\begin{defin}\rm\label{precube}
A \emph{precubical set} is a graded set $P = (P_n)_{n \geq 0}$ with  \emph{boundary operators} $d^k_i: P_n \to P_{n-1}$ $(n>0,\;k= 0,1,\; i = 1, \dots, n)$ satisfying the relations $d^k_i\circ d^l_{j}= d^l_{j-1}\circ d^k_i$ $(k,l = 0,1,\; i<j)$. If there exists a largest $n$ such that $P_n \not= \emptyset$, then this $n$ is called the \emph{dimension} of $P$. The degree of an element $x$ of $P$ will be denoted by $|x|$. The elements of degree $0$ are also called the \emph{vertices} of $P$. A morphism of precubical sets is a morphism of graded sets which is compatible with the boundary operators. The category of precubical sets will be denoted by $\Box \textbf{Set}$.
\end{defin} 

The category $\Box \textbf{Set}$ can be seen as the presheaf category of functors \linebreak $\Box^{op} \to \bf Set$ where $\Box$ is the small subcategory of $\bf Top$ whose objects are the standard $n$-cubes $I^n$ $(n \geq 0)$ and whose non-identity morphisms are composites of the maps $\delta^k_i\colon I^n\to I^{n+1}$ ($n \geq 0$, $i \in  \{1, \dots, n+1\}$, $k \in \{0,1\}$) given by $\delta_i^k(u_1,\dots, u_n)= (u_1,\dots, u_{i-1},k,u_i \dots, u_n)$. Here, we consider the $0$-cube as the one-point space $I^0 = \{()\}$. The \emph{precubical $n$-cube} is the $n$-dimensional precubical set $\mathbb I^n = \Box(-,I^n)$. By Yoneda's Lemma, an element $x$ of degree $n$ of a precubical set $P$ determines a unique morphism of precubical sets $x_{\sharp}\colon {\mathbb I}^n \to P$ such that $x_{\sharp}(id_{I^n}) = x$.

\begin{defin}\rm
Let $P$ be a precubical set and $x \in P_n$ be an element. We say that $x$ is \emph{regular} if the morphism $x_{\sharp}$ is injective. 
\end{defin}

We remark that if $x$ is regular, then so is each $d_i^kx$.

\begin{defin}\rm
A \emph{precubical subset} of a precubical set $P$ is a precubical set $Q$ such that $Q_n \subseteq P_n$ for all $n \geq 0$ and the boundary operators of $Q$ and $P$ coincide on $Q$. The \emph{opposite precubical set} of a precubical set $P$ with boundary operators $d_i^k$ is the precubical set $P^{op}$ with boundary operators $\partial_i^k$ defined by $P^{op}_r = P_r$ and $\partial_i^k = d_i^{1-k}$. The \emph{transposed precubical set} of a precubical set $P$ with boundary operators $d_i^k$ is the precubical set $P^t$ with boundary operators $\partial_i^k$ defined by $P^t_r = P_r$ and $\partial_i^k = d_{r+1-i}^k\colon P_r \to P_{r-1}$ ($i \in \{1, \dots, r\}$). The maps $P \mapsto P^{op}$ and $P \mapsto P^t$ extend to functors in the obvious way.
\end{defin}

The terminology of opposite and transposed precubical sets is an adaptation of the one used in \cite{GrandisBook} for cubical sets. Note that the opposite and transposed precubical set functors are involutions and that they preserve precubical subsets. Note also that a regular element of a precubical set $P$ is also regular as an element of $P^{op}$ and $P^t$.\\ 

Precubical sets can be realized geometrically as d-spaces in the sense of M. Grandis \cite{GrandisBook}. These spaces are defined as follows:

\begin{defin}\rm \cite[I.1.4]{GrandisBook}
A \emph{d-space} is a topological space $X$ together with a subset $dX \subseteq X^I$ such that
\begin{itemize}
\item[(i)] $dX$ contains all constant paths,
\item[(ii)] $dX$ is closed under composition with (not necessarily strictly) increasing maps $I \to I$,
\item[(iii)] $dX$ is closed under concatenation. 
\end{itemize}
The elements of $dX$ are called \emph{d-paths} in $X$. A \emph{d-map} is a continuous map $f\colon X \to Y$ between d-spaces such that for any $\omega \in dX$, $f\circ \omega \in dY$. The category of d-spaces and d-maps is denoted by d\textbf{Top}.
\end{defin}

\begin{ex}\rm
The \emph{directed interval} is the d-space $\vec I = (I,d\vec I)$ where $d\vec I$ consists of the (not necessarily strictly) increasing maps $I \to I$. The d-paths in a d-space $X$ are precisely the d-maps $\vec I \to X$. Taking the $n$-fold product of $\vec I$ with itself we obtain the \emph{directed $n$-cube} $\vec I^n$. Note that 
the product of two d-spaces $X$ and $Y$ is a d-space with respect to the set $d(X\times Y)$ corresponding to $dX\times dY$ under the bijection $(X\times Y)^I \approx X^I \times Y^I$.
\end{ex}

\begin{defin} \label{geomreal} \rm (cp. \cite{FajstrupGR}, \cite{GaucherGoubault}, \cite{Goubault}, \cite[I.1.6.7]{GrandisBook}, \cite{RaussenTraceCub}) The \emph{geometric realization} of a precubical set $P$ is the quotient space $|P|=(\coprod _{n \geq 0} P_n \times I^n)/\sim$ where the sets $P_n$ are considered as discrete spaces and the equivalence relation is given by
$$(d^k_ix,u) \sim (x,\delta_i^k(u)), \quad  x \in P_{n+1},\; u\in I^n,\; i \in  \{1, \dots, n+1\},\; k \in \{0,1\}.$$ 
The geometric realization $|P|$ is a d-space with respect to the set $d|P|$ consisting of increasing reparametrizations of finite concatenations of paths $\omega \colon I \to |P|$ of the form $\omega (t) = [x,\alpha(t)]$ where $x \in P_n$ and $\alpha$ is a continuous map $I \to I^n$ which is order-preserving with respect to the natural order of $I$ and the componentwise natural order of $I^n$. The geometric realization of a morphism of precubical sets $f\colon P \to Q$ is the d-map $|f|\colon |P| \to |Q|$ given by $|f|([x,u])= [f(x),u]$. With these definitions the geometric realization is a functor $|\,|\colon {\Box \bf Set} \to {\rm d}{\bf Top}$. 
\end{defin} 

\begin{ex}\rm
The map $\vec I^ n \to |\mathbb I^ n|$, $u \mapsto [id_{I^n},u]$ is an isomorphism of d-spaces. 
\end{ex}

We remark that the geometric realization of a precubical set $P$ is a CW-complex (cf. \cite{GaucherGoubault}). The $n$-skeleton of $|P|$ is the geometric realization of the $n$-dimensional precubical subset $P_{\leq n}$ of $P$ defined by $(P_{\leq n})_m = P_m$ $(m\leq n)$. The closed $n$-cells of $|P|$ are the d-spaces $|x_{\sharp}(\mathbb I^n)|$ where $x \in P_n$. The characteristic map of the cell $|x_{\sharp}(\mathbb I^n)|$ is the d-map $\vec I^n \stackrel{\cong}{\to} |\mathbb I ^n| \stackrel{|x_{\sharp}|}{\to} |P|$ and this map is an isomorphism onto its image if and only if $x$ is regular. Note also that a natural isomorphism of d-spaces $\sigma_P \colon |P| \to |P^t|$ is given by $[z,(u_1,\dots,u_r)] \mapsto [z,(u_r, \dots, u_1)]$. Note finally that the geometric realization of a precubical subset $Q$ of $P$ is both a CW-subcomplex of $|P|$ and a d-subspace of $|P|$ in the sense of the following definition:

\begin{defin}\rm \cite[I.1.4.1]{GrandisBook}
A \emph{d-subspace} of a d-space $X$ is a d-space $A$ such that the topological space $A$ is a subspace of $X$ and $dA = \{\omega \in A^I\,|\, (A\hookrightarrow X) \circ \omega \in dX\}$.
\end{defin}

Reversing the direction of the d-paths of a d-space one obtains the opposite d-space:

\begin{defin}\rm \cite[I.1.4.0]{GrandisBook} 
Given a path $\omega \colon I \to X$ we denote the inverse path $I \to X$, $t \mapsto \omega(1-t)$ by $\bar \omega$. The \emph{opposite d-space} of a d-space $X = (X,dX)$ is the d-space $X^{op} = (X,dX^{op})$ defined by $dX^{op} = \{\bar \omega \,|\, \omega \in dX\}$. The map $X \mapsto X^{op}$ extends to a functor in the obvious way.
\end{defin}

Note that the opposite d-space functor is an involution. Note also that if $A$ is a d-subspace of $X$, then $A^ {op}$ is a d-subspace of $X^{op}$. Note finally that for a precubical set $P$, a natural isomorphism of d-spaces $\phi_P \colon |P|^{op} \to |P^{op}|$ is given by $[z,(u_1,\dots,u_r)] \mapsto [z,(1-u_1, \dots, 1-u_r)]$.

\section{Dihomotopy, the fundamental bipartite graph, and $\rm d$-path spaces}

We shall work with the following notion of directed homotopy:

\begin{defin} \label{defdihomotopy} \rm
Two d-maps $f,g\colon X\to Y$ are said to be \emph{dihomotopic} if there exists a homotopy $H\colon X \times I \to Y$ from $f$ to $g$ such that each map $H(-,t)$ is a d-map. Such a homotopy is called a \emph{dihomotopy} from $f$ to $g$. If $f$ and $g$ coincide on a d-subspace $A \subseteq X$, then $f$ and $g$ are said to be \emph{dihomotopic relative to $A$} if there exists a \emph{dihomotopy relative to $A$} from $f$ to $g$, i.e., a  dihomotopy $H\colon X \times I \to Y$ from $f$ to $g$ such that each map $H(-,t)$ coincides with $f$ and $g$ on $A$. Let $X$ and $Y$ be d-spaces with a common d-subspace $A$. A d-map $f\colon X \to Y$ satisfying $f(a) = a$ for all $a \in A$ is said to be a \emph{dihomotopy equivalence relative to $A$} if there exists a \emph{dihomotopy inverse relative to $A$} of $f$, i.e., a d-map $g\colon Y \to X$ such that $g(a) = a$ for all $a \in A$ and such that $g\circ f$ and $f\circ g$ are dihomotopic relative to $A$ to the identities of $X$ and $Y$, respectively. A \emph{dihomotopy equivalence} is a d-map which is a  dihomotopy equivalence relative to the empty \linebreak d-space. Two d-spaces $X$ and $Y$ with a common d-subspace $A$ are said to be \emph{dihomotopy equivalent relative to $A$} if there exists a dihomotopy equivalence relative to $A$ between them. Two d-spaces $X$ and $Y$ are \emph{dihomotopy equivalent} if they are dihomotopy equivalent relative to the empty d-space.
\end{defin}

We remark that (relative) dihomotopy is an equivalence relation which is compatible with the composition of d-maps.  Some authors, as for instance M. Grandis \cite{GrandisBook}, work with a stronger notion of directed homotopy, called \emph{d-homotopy}, where the homotopies are required to be d-maps $X \times \vec I \to Y$. The reader is referred to L. Fajstrup \cite{Fajstrup} for a result concerning the equivalence of the two notions of directed homotopy for directed paths.\\

The one-dimensional information of a d-space is contained in its fundamental category, which is the directed analogue of the fundamental groupoid of a topological space.

\begin{defin}\rm (\cite{Goubault}, \cite[I.3]{GrandisBook})
The \emph{fundamental category} of a d-space $X$ is the category $\vec \pi_1(X)$ defined as follows: The objects are the elements of $X$ and the set of morphisms from an element $x$ to an element $y$ is the set of dihomotopy classes relative to $\{0,1\}$ of d-paths from $x$ to $y$. The map $X \mapsto \vec \pi_1(X)$ extends in the obvious way to a functor from d$\bf Top$ to the category of small categories.
\end{defin}

The fundamental category of a d-space is obviously a huge object. This led to the development of several methods to extract the essential directed information of the fundamental category (cf. \cite{BubenikExtremal}, \cite{FajstrupGR}, \cite{Components}, \cite{Components2}, \cite[I.3]{GrandisBook}, \cite{RaussenInvariants}). In \cite{BubenikExtremal}, P. Bubenik introduced the fundamental bipartite graph of a d-space:

\begin{defin}\rm \cite{BubenikExtremal}
An element $a$ of a d-space $X$ is said to be \emph{minimal} (\emph{maximal}) if any morphism in $\vec \pi_1(X)$ with target (source) $a$ has source (target) $a$. An element of $X$ is \emph{extremal} if it is minimal or maximal. The d-subspace of $X$ consisting of the extremal elements is denoted by $Extrl(X)$. The \emph{fundamental bipartite graph} of a d-space $X$, denoted by $\vec \pi_1(X,Extrl(X))$, is the full subcategory of $\vec \pi_1(X)$ generated by $Extrl(X)$. 
\end{defin}

Note that the fundamental bipartite graph of a d-space is a bipartite graph if one ignores the identity morphisms. In a d-space model of a concurrent system, initial states of the system are modeled by minimal elements and final states and deadlocks are modeled by maximal elements. The fundamental bipartite graph of the d-space represents the essential execution schedules between these critical states of the system. For precubical sets there is another definition of minimal, maximal, and extremal elements:

\begin{defin}\rm
Let $P$ be a precubical set. An element $v\in P_0$ is said to be \emph{minimal} (\emph{maximal}) if there is no element $x \in P_1$ such that $d_1^1x = v$ ($d_1^0x = v$). An element of $P$ is \emph{extremal} if it is minimal or maximal. The $0$-dimensional precubical subset of $P$ consisting of the extremal elements is denoted by $Extrl(P)$.  
\end{defin}

We remark that $Extrl(P^{op}) = Extrl(P^t) = Extrl(P)$. The easy proof of the following proposition is left to the reader.

\begin{prop} \label{extrlP}
Let $P$ be a precubical set. An element $a \in |P|$ is minimal (maximal) if and only if there exists a minimal (maximal) element $v\in P_0$ such that $a = [v,()]$. Consequently, $|Extrl(P)| = Extrl(|P|)$.
\end{prop}

In \cite{FahrenbergRaussen} and \cite{RaussenInvariants}, U. Fahrenberg and M. Rau\ss en introduce \emph{trace spaces}, which are quotient spaces of the d-path spaces we define next. In \cite{RaussenTraceCub}, it is shown that trace and d-path spaces are actually homotopy equivalent for geometric realizations of precubical sets.

\begin{defin} \rm
Let $X$ be a d-space and $a, b \in X$. The \emph{d-path space} $\vec {\mathcal P} (X)(a,b)$ is defined to be the subspace of $X^I$ consisting of the d-paths from $a$ to $b$.
\end{defin}

\begin{rem} \label{pathcomp}\rm Note that the path components of $\vec {\mathcal P} (X)(a,b)$ are the morphisms from $a$ to $b$ in $\vec \pi_1 (X)$.
\end{rem}

We have the following result on the dihomotopy invariance of the fundamental bipartite graph and d-path spaces:

\begin{theor} \label{FBGtheor}
Let $X$ and $Y$ be two d-spaces such that $Extrl(X) = Extrl(Y)$ and let $f\colon X \to Y$ be a dihomotopy equivalence relative to $Extrl(X)$. Then the functor $\vec \pi_1(f)\colon \vec \pi_1(X) \to \vec \pi_1(Y)$ restricts to an isomorphism of fundamental bipartite graphs $\vec \pi_1(X,Extrl(X)) \to \vec \pi_1(Y,Extrl(Y))$. Moreover, for any two elements $a, b \in Extrl(X)$ the map $f^I\colon X^I \to Y^I$, $\omega \mapsto f\circ \omega$ restricts to a homotopy equivalence $\vec {\mathcal P}(X)(a,b) \to \vec {\mathcal P}(Y)(a,b)$.
\end{theor}

\dem
By \ref{pathcomp}, the statement on the fundamental bipartite graphs follows from the statement on d-path spaces. In order to prove the latter, let $a, b \in Extrl(X) = Extrl(Y)$ and $g\colon Y \to X$ be a dihomotopy inverse relative to $Extrl(X)$ of $f$. Let $f_*\colon \vec {\mathcal P}(X)(a,b) \to \vec {\mathcal P}(Y)(a,b)$ and $g_*\colon \vec {\mathcal P}(Y)(a,b) \to \vec {\mathcal P}(X)(a,b)$ be the restrictions of $f^I$ and $g^I$. Let $H\colon X \times I \to X$ be a dihomotopy relative to $Extrl(X)$ from $id_X$ to $g\circ f$. Then a homotopy $h\colon \vec {\mathcal P}(X)(a,b) \times I \to \vec {\mathcal P}(X)(a,b)$ from $id_{\vec {\mathcal P}(X)(a,b)}$ to $g_*\circ f_*$ is given by $h(\omega ,t)(s) = H(\omega (s),t)$. Similarly, $id_{\vec {\mathcal P}(Y)(a,b)} \simeq f_*\circ g_*$.
\findem

\begin{cor} \label{FBG}
Let $P$ and $Q$ be two precubical sets such that $Extrl(P) = Extrl(Q)$ and let $f\colon |P| \to |Q|$ be a dihomotopy equivalence relative to $|Extrl(P)|$. Then the functor $\vec \pi_1(f)\colon \vec \pi_1(|P|) \to \vec \pi_1(|Q|)$ restricts to an isomorphism of fundamental bipartite graphs $\vec \pi_1(|P|,Extrl(|P|)) \to \vec \pi_1(|Q|,Extrl(|Q|))$. Moreover, for any two elements $v, w \in Extrl(P)$ the map $f^I\colon |P|^I \to |Q|^I$, $\omega \mapsto f\circ \omega$ restricts to a homotopy equivalence $\vec {\mathcal P}(|P|)([v,()],[w,()]) \to \vec {\mathcal P}(|Q|)([v,()],[w,()])$.
\end{cor}

The last proposition of this section permits us to dualize results on dihomotopy equivalences for precubical sets. The straightforward proof is left to the reader.

\begin{prop} \label{dualization}
Let $P$ and $Q$ be two precubical sets with a common precubical subset $R$. \begin{itemize}
\item[(i)] If $f\colon |P^{op}| \to |Q^{op}|$ is a dihomotopy equivalence relative to $|R^{op}|$, then \linebreak  $(\phi_Q^{-1}\circ f \circ \phi_P)^{op}\colon |P| \to |Q|$ is a dihomotopy equivalence relative to $|R|$.
\item[(ii)] If $f\colon |P^t| \to |Q^t|$ is a dihomotopy equivalence relative to $|R^t|$, then \linebreak $\sigma_Q^{-1}\circ f \circ \sigma_P\colon |P| \to |Q|$ is a dihomotopy equivalence relative to $|R|$. 
\end{itemize}
\end{prop}

\section{Construction of $\rm d$-maps}

The purpose of this section is to present a method to construct d-maps and dihomotopies between geometric realizations of precubical sets. The maps and homotopies we consider in the next sections can be shown to be d-maps and dihomotopies by checking the conditions we establish in this section.

\begin{defin}\rm
A subset $Z$ of a partially ordered set $(J,\leq)$ is called \emph{order-convex} if for any two elements $a,b \in Z$, $\{z \in J\,| a\leq z \leq b\} \subseteq Z$.
\end{defin}

\begin{rem}\rm \label{ordercon}
If $\alpha \colon I \to I^ m$ is an order-preserving map and $s \leq t$ are elements of $I$ such that $\alpha(s)$ and $\alpha(t)$ belong to an order-convex set $Z \subseteq I^m$, then $[s,t] \subseteq \alpha^{-1}(Z)$.
\end{rem}

\begin{prop} \label{d-maps}
Let $P$ and $Q$ be precubical sets and $f\colon \coprod \limits_{r \geq 0} P_r \times I^r \to |Q|$ be a continuous map. Suppose that
\begin{itemize}
\item[(i)] $f(d_i^kz,u) = f(z, \delta_i^ku)$ for all $r\geq 1$, $z \in P_r$, $u \in I^{r-1}$, $i \in \{1,\dots r\}$, $k \in \{0,1\}$,
\item[(ii)] for all $r \geq 1$ and $z \in P_r$ there exist a finite closed order-convex covering $\mathcal A_z$ of $I^r$, a function $\zeta_z \colon \mathcal A_z \to \coprod \limits_{m \geq 0} Q_m$, and a family of order-preserving maps $\{f_{z,Z}\colon Z \to I^{|\zeta_z(Z)|}\}_{Z \in \mathcal A_z}$ such that for all $Z \in \mathcal A_{z}$ and $u \in Z$, $f(z,u) = [\zeta_z(Z),f_{z,Z}(u)]$. 
\end{itemize}
Then a d-map $\bar f \colon |P| \to |Q|$ is given by $\bar f([z,u]) = f(z,u)$.
\end{prop}

\dem
By condition (i), $\bar f$ is well-defined and continuous. Let $r \geq 0$ and $z \in P_r$ and consider a path $\omega \in d|P|$ of the form $\omega(t) = [z,\alpha(t)]$ where the map $\alpha\colon I \to I^{r}$ is order-preserving. We have to show that $\bar f \circ \omega \in d|Q|$. If $r= 0$, then $\bar f \circ \omega$ is a constant path and therefore $\bar f \circ \omega \in d|Q|$. Let  $r \geq 1$. Let $\mathcal B_z$ be the subset of $\mathcal A_z$ consisting of the sets $Z \in \mathcal A_z$ such that $\alpha^{-1}(Z)$ has more than one element. Then $I = \bigcup \limits_{Z \in \mathcal B_z} \alpha^{-1}(Z)$. Indeed, else there would exist an element $s \in I \setminus \bigcup \limits_{Z \in \mathcal B_z} \alpha^{-1}(Z)$ and one would have $\bigcup \limits_{\substack{Z \in \mathcal A_z\\s\in \alpha^{-1}(Z)}} \alpha^{-1}(Z) = \{s\}$ and hence $I\setminus \{s\} = \bigcup \limits_{\substack{Z \in \mathcal A_z\\s\notin \alpha^{-1}(Z)}}\alpha^{-1}(Z)$ which is impossible since $I \setminus \{s\}$ is not closed in $I$. Define a subset $\{Z_1, \dots, Z_l\} \subseteq \mathcal B_z$ such that $0 < \max \alpha^{-1}(Z_1) < \cdots < \max \alpha^{-1}(Z_l) = 1$, $[0,\max \alpha^{-1}(Z_1)] = \alpha^{-1}(Z_1)$, and $[\max \alpha^{-1}(Z_{i-1}),\max \alpha^{-1}(Z_i)] \subseteq \alpha^{-1}(Z_i)$ for $i \in \{2,\dots, , l\}$ inductively as follows. Choose $Z_1 \in \mathcal B_z$ such that $0 \in \alpha^{-1}(Z_1)$. Then $0 < \max \alpha^{-1}(Z_1)$ and, by \ref{ordercon}, $[0,\max \alpha^{-1}(Z_1)] = \alpha^{-1}(Z_1)$. Suppose that $Z_{i}$ has been defined and that $\max \alpha^{-1}(Z_{i}) < 1$. Then $$[0,\max \alpha^{-1}(Z_{i})] = \bigcup \limits_{\substack{Z \in \mathcal B_z\\\max \alpha^{-1}(Z)\leq \max \alpha^{-1}(Z_i)}}\alpha^{-1}(Z)$$ and hence $$ ]\max \alpha^{-1}(Z_i),1] \subseteq \bigcup \limits_{\substack{Z \in \mathcal B_z\\\max \alpha^{-1}(Z) > \max \alpha^{-1}(Z_i)}}\alpha^{-1}(Z).$$
Since this is a finite union of closed subsets of $I$, it even contains the closed interval $[\max \alpha^{-1}(Z_i),1]$ as a subset. Therefore we may choose  $Z_{i+1} \in \mathcal B_z$ such that $\max \alpha^{-1}(Z_i) < \max \alpha^{-1}(Z_{i+1})$ and $\max \alpha^{-1}(Z_i) \in \alpha^{-1}(Z_{i+1})$. By \ref{ordercon}, we then also have $[\max \alpha^{-1}(Z_i),\max \alpha^{-1}(Z_{i+1})] \subseteq \alpha^{-1}(Z_{i+1})$. Since $\mathcal B_z$ is finite, the process terminates after a finite number of steps. Set $b_i = \max \alpha^{-1}(Z_i)$ $(i = 1, \dots, l)$ and $b_0 = 0$. It suffices to show that for each $i \in \{1,\dots ,l\}$, the path $\gamma_i\colon I \to |Q|$, $t \mapsto \bar f \circ \omega((1-t)b_{i-1} + tb_i)$ belongs to $d|Q|$. Let $\beta_i$ be the composite $$I \stackrel{(1-t)b_{i-1}+tb_i}{\longrightarrow} [b_{i-1},b_i] \stackrel{\alpha}{\to}Z_i \stackrel{f_{z,Z_i}}{\to} I^{|\zeta_z(Z_i)|}.$$ Then $\beta_i$ is order-preserving and $\gamma_i (t) = [\zeta_z(Z_i),\beta_i(t)]$ $(t\in I)$. Thus, $\gamma_i \in d|Q|$. 
\findem

\begin{prop} \label{dihom}
Let $P$ and $Q$ be precubical sets and $h\colon \coprod \limits_{r \geq 0} P_r \times I^r \times I\to |Q|$ be a continuous map. Suppose that
\begin{itemize}
\item[(i)] $h(d_i^kz,u,t) = h(z, \delta_i^ku,t)$ for all $r\geq 1$, $z \in P_r$, $u \in I^{r-1}$, $i \in \{1,\dots r\}$, $k \in \{0,1\}$, $t\in I$,
\item[(ii)] for all $t\in I$, $r \geq 1$, and $z \in P_r$ there exist a finite closed order-convex covering $\mathcal A_{z,t}$ of $I^r$, a function $\zeta_{z,t} \colon \mathcal A_{z,t} \to \coprod \limits_{m \geq 0} Q_m$, and a family of order-preserving maps $\{h_{z,t,Z}\colon Z \to I^{|\zeta_{z,t}(Z)|}\}_{Z \in \mathcal A_{z,t}}$ such that for all $Z \in \mathcal A_{z,t}$ and $u \in Z$, $h(z,u,t) = [\zeta_{z,t}(Z),h_{z,t,Z}(u)]$. 
\end{itemize}
Then a dihomotopy $H \colon |P|\times I \to |Q|$ is given by $H([z,u],t) = h(z,u,t)$.
\end{prop}

\dem
By condition (i), $H$ is well-defined and continuous. By \ref{d-maps}, the map \linebreak $H(\;,t)\colon |P| \to |Q|$ is a d-map for each $t \in I$.
\findem

\section{One-dimensional reduction}

In ordinary homotopy theory, a contractible subspace of a topological space can be collapsed to a point, at least if the space and the subspace form an NDR-pair. The resulting quotient space has the same homotopy type as the original space. In directed homotopy theory, the situation is more complicated. Consider, for example, the geometric realization of a precubical set $P$ with two vertices and two edges that looks like the graph in figure \ref{fig1b}. If one collapses one of the edges to a point, one obtains the directed circle $\vec S^1$, which is the geometric realization of a precubical set with one vertex and one edge. The d-spaces $|P|$ and $\vec S^1$ are \emph{not} dihomotopy equivalent. The following theorem gives a condition under which it is possible to collapse an edge in the geometric realization of a precubical set to a point without changing the directed homotopy type relative to the extremal elements: 

\begin{theor} \label{elimination3}
Let $P$ be a precubical set, $b \in \{0,1\}$, and $x \in P_1$ be a regular element such that
\begin{itemize}
\item[(i)] there is no element $y \in P_1 \setminus \{x\}$ such that $d_1^{1-b}y = d_1^{1-b}x$,
\item[(ii)] no element in $P_1$ having $d_1^{1-b}x$ in its boundary  belongs to the boundary of an element in $P_2$. 
\end{itemize}
Consider the set $Y = \{y \in P_1\,|\, d_1^by = d_1^{1-b}x\}$. Then a precubical set $Q$ such that $Y \subseteq Q_1$, $Q \setminus Y$ is a common precubical subset of $P$ and $Q$, and $|P|$ and $|Q|$ are dihomotopy equivalent relative to $|Q \setminus Y|$ is given by $Q_0 = P_0 \setminus \{d_1^{1-b}x\}$, $Q_1 = P_1 \setminus \{x\}$, $Q_r = P_r$ $(r > 1)$, and the boundary operators $D _i^k$ defined by 
$$ D _i^kz= \left\{ \begin{array}{ll}
d_1^bx, & z\in Y,\, i=1,\,k=b,\\
d_i^kz, & \mbox{else}.
\end{array}\right.$$
Moreover, if $Y \not= \emptyset$, then $Extrl(P) = Extrl(Q) \subseteq Q\setminus Y$ and $|P|$ and $|Q|$ have isomorphic fundamental bipartite graphs and homotopy equivalent d-path spaces for each pair of extremal elements.
\end{theor}

\begin{figure}
\begin{center}
\subfigure[$P$]
{ 
\includegraphics[height=2.5cm]{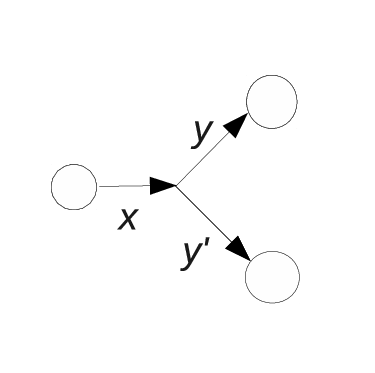}
\label{thm1P}
}
\hspace{2cm}
\subfigure[$Q$]
{
\includegraphics[height=2.5cm]{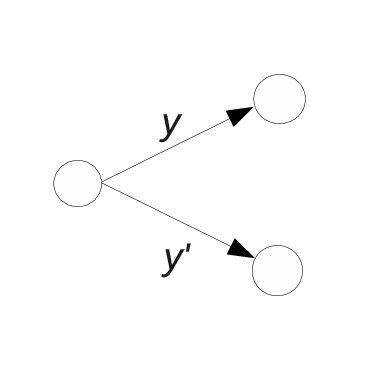}
\label{thm1Q}
}
\caption{Precubical sets $P$ and $Q$ in theorem \ref{elimination3}} \label{thm1}
\end{center}
\end{figure}

\dem
We only consider the case $b=0$. The case $b=1$ is dual and can be deduced from the case $b=0$ using \ref{dualization}(i). Note first that the regularity of $x$ and the two conditions of the theorem ensure that $Q$ is a precubical set, $Y \subseteq Q_1$, and $Q \setminus Y$ is a common precubical subset of $P$ and $Q$. The situation is illustrated in figure \ref{thm1}, where $Y = \{y,y'\}$ and $Q\setminus Y$ is represented by the circles. 

Consider the continuous maps $f\colon \coprod \limits_{r\geq 0} P_r\times I^r \to |Q|$ and $g\colon \coprod \limits_{r\geq 0} Q_r\times I^r \to |P|$ defined by $$f(z,u) = \left\{ \begin{array}{ll}
[d_1^0x,()], & z \in \{x,d_1^1x\},\\   

[z,u], & z \notin \{x,d_1^1x\}
\end{array}\right.$$
and 
$$g(z,u) = \left\{ \begin{array}{ll}
[x,2u], & z \in Y,\,u\leq 1/2\\   

[z,2u-1], & z \in Y,\,u\geq 1/2\\

[z,u], & z \notin Y.
\end{array}\right.$$
It is straightforward to check that $f$ and $g$ satisfy the conditions of \ref{d-maps}. It follows that the maps $\bar f \colon |P| \to |Q|$, $\bar f([z,u]) = f(z,u)$ and $\bar g \colon |Q| \to |P|$, $\bar g([z,u]) = g(z,u)$ are d-maps. Note that both $\bar f$ and $\bar g$ restrict to the identity on $|Q \setminus Y|$. 

We show that $\bar f$ and $\bar g$ are inverse dihomotopy equivalences relative to $|Q \setminus Y|$. Consider the map $\phi\colon \coprod \limits_{r\geq 0} P_r\times I^r \times I \to |P|$ defined by $$\phi(z,u,t) = \left\{ \begin{array}{ll}
[z,u], & z \in\{x,d_1^1x\},\, t\leq 1/2\\   

[x,(2-2t)u], & z = x,\, t\geq 1/2,\\

[x,2-2t], & z = d_1^1x,\, t\geq 1/2,\\

[z,(1-2t)u], & z \in Y,\,u\leq 1/2,\, t\leq 1/2,\\

[x,2-2t+(2t-1)2u], & z \in Y,\,u\leq 1/2,\, t\geq 1/2,\\

[z,(1-2t)u+2t(2u-1)], & z \in Y,\,u\geq 1/2,\, t\leq 1/2,\\

[z,2u-1], & z \in Y,\,u\geq 1/2,\, t\geq 1/2,\\

[z,u], & z \notin Y\cup\{x,d_1^1x\}.
\end{array}\right.$$
It is straightforward to check that $\phi$ is well-defined and continuous and that it satisfies the conditions of \ref{dihom}. Therefore the map $\Phi \colon |P| \times I \to |P|$, $\Phi([z,u],t) = \phi(z,u,t)$ is a dihomotopy. We have $\Phi([z,u],0) = [z,u]$ and $\Phi ([z,u],1) = \bar g\circ \bar f ([z,u])$. Moreover, $\Phi([z,u],t) = [z,u]$ for all $[z,u] \in |P \setminus (Y\cup \{x,d_1^1x\})| = |Q\setminus Y|$ and $t\in I$. It follows that $\Phi$ is a dihomotopy relative to $|Q\setminus Y|$ from $id_{|P|}$ to $\bar g \circ \bar f$.

Consider the continuous map $\psi\colon \coprod \limits_{r\geq 0} Q_r\times I^r \times I \to |Q|$ defined by $$\psi(z,u,t) = \left\{ \begin{array}{ll}
[z,(1-t)u], & z \in Y,\,u\leq 1/2,\\

[z,(1-t)u+t(2u-1)], & z \in Y,\,u\geq 1/2,\\

[z,u], & z \notin Y.
\end{array}\right.$$
One easily verifies that $\psi$ satisfies the conditions of \ref{dihom}. Therefore the map \linebreak $\Psi \colon |Q| \times I \to |Q|$, $\Psi([z,u],t) =\psi(z,u,t)$ is a dihomotopy. We have $\Psi([z,u],0) = [z,u]$ and $\Psi ([z,u],1) = \bar f\circ \bar g ([z,u])$. Moreover, $\Psi([z,u],t) = [z,u]$ for all $[z,u] \in |Q \setminus Y|$ and $t\in I$. It follows that $\Psi$ is a dihomotopy relative to $|Q\setminus Y|$ from $id_{|Q|}$ to $\bar f \circ \bar g$.

One easily shows that $Extrl(P) = Extrl(Q) \subseteq  Q\setminus Y$ if $Y \not= \emptyset$. By \ref{FBG}, this implies that $|P|$ and $|Q|$ have isomorphic fundamental bipartite graphs and homotopy equivalent d-path spaces for each pair of extremal elements. 
\findem

\begin{rem}\rm
Note that if the set $Y$ has exactly one element, the equivalence between $|P|$ and $|Q|$ can be seen as a T-homotopy equivalence in the sense of \cite{GaucherGoubault}.
\end{rem}

\section{Two-dimensional reduction}

In opposition to the situation in ordinary homotopy theory, the removal of a cube and a free face from a precubical set changes in general the directed homotopy type of the geometric realization. In this section we prove two theorems which give conditions under which it is possible to eliminate a 2-dimensional cube and one or two free faces in a precubical set without changing the directed homotopy type relative to the extremal elements of the geometric realization. 

\begin{theor}\label{elimination2}
Let $P$ be a precubical set, $a\in\{1,2\}$, $b \in \{0,1\}$, and $x \in P_2$ be a regular element such that 
\begin{itemize}
\item[(i)] no element of $P_2\setminus \{x\}$ has $d_a^{1-b}x$ or $d_{3-a}^bx$ in its boundary, 
\item[(ii)] there is no element $y \in P_1\setminus \{d_a^{1-b}x\}$ such that $d_1^by = d_1^bd_a^{1-b}x$,
\item[(iii)] no element of the set $Y = \{y \in P_1\setminus \{d_{3-a}^bx\} \,|\, d_1^{1-b}y = d_1^{1-b}d_{3-a}^bx\}$ is in the boundary of some element in $P_2$.
\end{itemize}
Then a precubical subset $Q$ of $P$ and a precubical subset $R$ of $Q$ such that the inclusion $\iota\colon |Q| \hookrightarrow |P|$ is a dihomotopy equivalence relative to $|R|$ are given by $Q_0 = P_0$, $Q_1 = P_1\setminus \{d_{3-a}^bx\}$, $Q_2 = P_2\setminus \{x\}$, $Q_r = P_r$ $(r >2)$, $R_0 = Q_0 \setminus \{d_1^{1-b}d_{3-a}^bx\}$, $R_1= Q_1 \setminus (\{d_a^{1-b}x\}\cup Y)$, and $R_r= Q_r$ $(r\geq 2)$. Moreover, if $Y \not= \emptyset$, then $Extrl(P) = Extrl(Q) \subseteq R$ and $\iota$ induces an isomorphism of fundamental bipartite graphs and a homotopy equivalence of d-path spaces for each pair of extremal elements.
\end{theor}

\begin{figure}
\begin{center}
\subfigure[$P$]
{ 
\includegraphics[height=3.5cm]{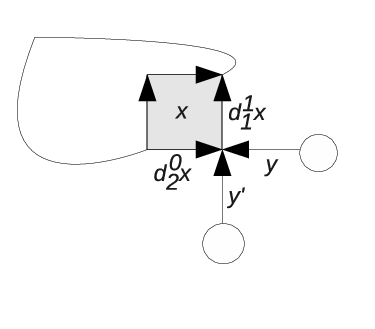}
\label{thm2P}
}
\hspace{2cm}
\subfigure[$Q$]
{
\includegraphics[height=3.5cm]{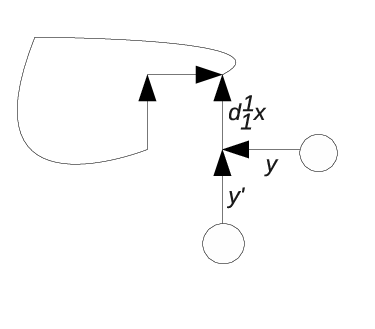}
\label{thm2Q}
}
\caption{Precubical sets $P$ and $Q$ in theorem \ref{elimination2}} \label{thm2}
\end{center}
\end{figure}

\dem
We only consider the case $a=1$ and $b= 0$. The remaining cases are dual and can be deduced from the case $a=1$ and $b= 0$ using \ref{dualization}. Note first that the regularity of $x$ and the three conditions of the theorem ensure that $Q$ is a precubical subset of $P$ and $R$ is a precubical subset of $Q$. The situation is illustrated in figure \ref{thm2}, where $Y = \{y,y'\}$ and $R$ is represented by the circles and the area limited by the curves.

\begin{figure}
\begin{center}
\includegraphics[height=5cm]{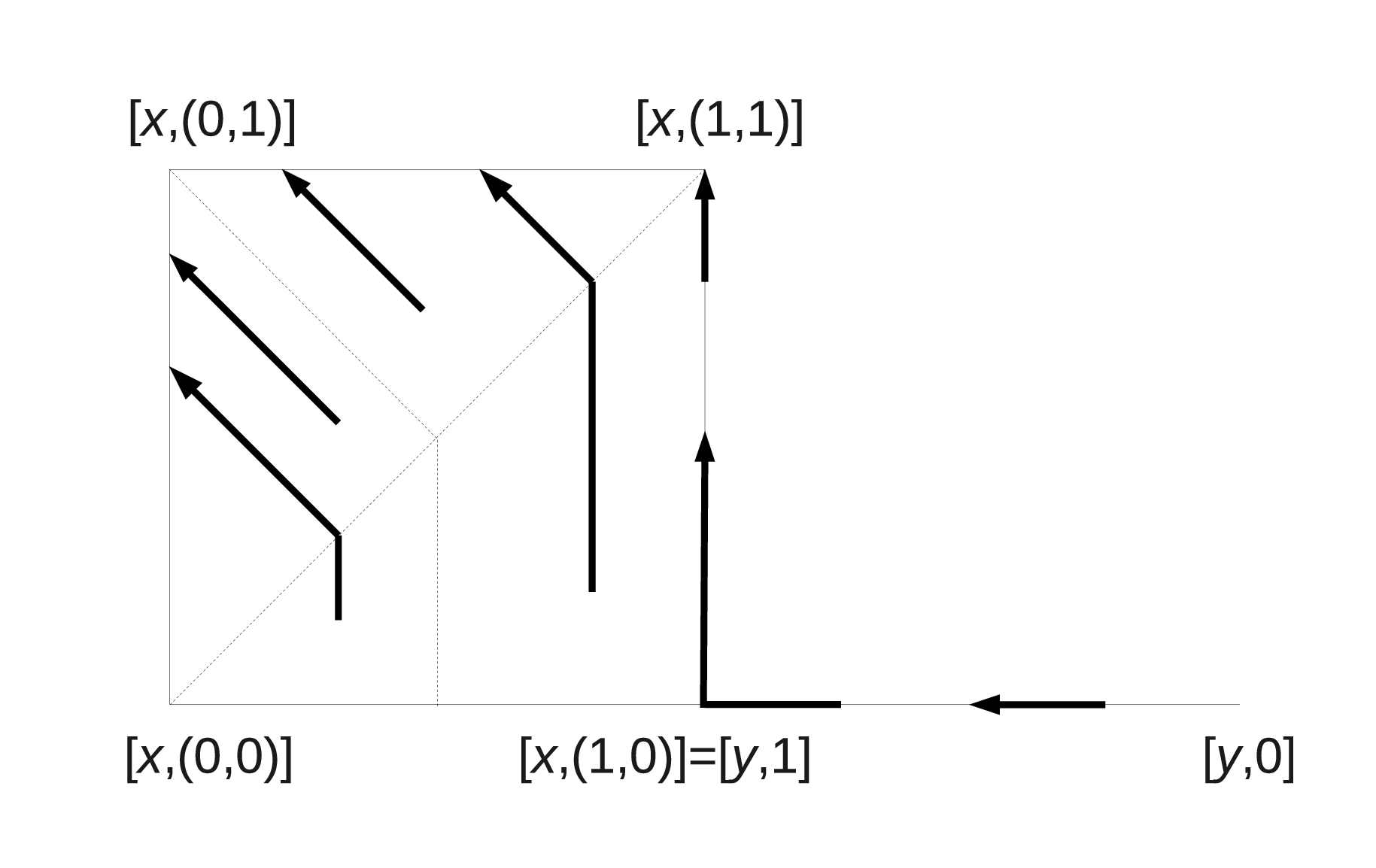}
\caption{The dihomotopy $H$ moves the beginning points of the thick lines to the endpoints of the arrows} \label{H}
\end{center}
\end{figure}

We construct a dihomotopy $H\colon |P| \times I \to |P|$ using \ref{dihom}. Consider the map $h\colon \coprod \limits_{r \geq 0} P_r \times I^r \times I\to |P|$ defined by \newline$h(x,(u_1,u_2),t)$
$$ = \left\{ \begin{array}{ll}
\vspace{1mm}
[x,(u_1,u_2)], & t \leq \frac{1}{3},\\

\vspace{1mm}
[x,(u_1,u_2)], & \frac{1}{3} \leq t \leq \frac{2}{3},\, u_1 \leq u_2,\\

\vspace{1mm}
[x,(u_1, (3t-1)u_1 + (2-3t)u_2)], & \frac{1}{3} \leq t \leq \frac{2}{3},\, u_1 \geq u_2,\\

\vspace{1mm}
[x,((3-3t)u_1,u_2+ (3t-2)u_1)], & t \geq \frac{2}{3},\, u_1 \leq u_2,\\& u_2 \leq 1 - u_1,\\

\vspace{1mm}
[x,(u_1+(3t-2)(u_2-1),(3-3t)u_2 +3t-2)], &  t \geq \frac{2}{3},\, u_1 \leq u_2,\\& u_2 \geq 1-u_1,\\

\vspace{1mm}
[x,((3-3t)u_1,(3t-1)u_1)], &  t \geq \frac{2}{3},\, u_1 \geq u_2,\, u_1 \leq \frac{1}{2},\\

[x,((3t-1)(u_1-1) +1,(3-3t)u_1 +3t-2)], &  t \geq \frac{2}{3},\, u_1 \geq u_2,\, u_1 \geq \frac{1}{2},
\end{array}\right.$$ 
$h(d_1^1x,u,t) = h(x,(1,u),t)$, $h(d_2^0x,u,t) = h(x,(u,0),t)$, $h(d_1^1d_2^0x,(),t) = \linebreak h(x,(1,0),t)$,
$$h(y,u,t) = \left\{ \begin{array}{ll}
\vspace{1mm}
[y,(1+3t)u], & t \leq \frac{1}{3},\, u \leq \frac{1}{2},\\

\vspace{1mm}
[y, (1-3t)u+3t], & t \leq \frac{1}{3},\, u \geq \frac{1}{2},\\

\vspace{1mm}
[y,2u], & t \geq \frac{1}{3},\, u \leq \frac{1}{2},\\

\vspace{1mm}
[d_1^1x,(3t-1)(2u-1)], & \frac{1}{3} \leq t \leq \frac{2}{3},\, u \geq \frac{1}{2},\\

[d_1^1x,2u-1], & t \geq \frac{2}{3},\, u \geq \frac{1}{2}\\
\end{array}\right.$$ 
for $y \in Y$, and $h(z,(u_1,\dots,u_r),t) = [z,(u_1,\dots,u_r)]$ for $z \in R$. It is straightforward to check that $h$ is well-defined and continuous. A tedious but also straightforward verification shows that the conditions of \ref{dihom} hold. Therefore a dihomotopy $H\colon |P| \times I \to |P|$ is given by $H([z,(u_1, \dots, u_r)],t) = h(z,(u_1,\dots,u_r),t)$. A picture indicating what the dihomotopy $H$ does is given in figure \ref{H}.
By definition, $H([z,(u_1, \dots, u_r)],t) = [z,(u_1, \dots, u_r)]$ for all $r\geq 0$, $z \in R_r$, $(u_1, \dots, u_r) \in I^r$, and $t \in I$. One easily checks that $H([z,(u_1, \dots, u_r)],0) = [z,(u_1, \dots, u_r)]$ and $H([z,(u_1, \dots, u_r)],1) \in |Q|$ for all $r \geq 0$, $z \in P_r$, and $(u_1, \dots, u_r) \in I^r$. One also checks easily that $H([z,(u_1,\dots,u_r)],t) \in |Q|$ for all $r \geq 0$, $z \in Q_r$, $(u_1, \dots, u_r) \in I^r$, and $t \in I$. This implies that the inclusion $\iota\colon |Q| \hookrightarrow |P|$ is a dihomotopy equivalence relative to $|R|$. Indeed, consider the d-map $f\colon |P| \to |Q|$ given by  $f([z,(u_1,\dots,u_r)]) = H([z,(u_1,\dots,u_r)],1)$. A dihomotopy relative to $|R|$ from $id_{|Q|}$ to $f\circ \iota$ is given by $G\colon |Q|\times I \to |Q|$, $G([z,(u_1,\dots,u_r)],t) = H([z,(u_1,\dots,u_r)],t)$ and $H$ is a dihomotopy relative to $|R|$ from $id_{|P|}$ to $\iota\circ f$.

Suppose now that $Y \not = \emptyset$. One easily sees that $Extrl(Q) = Extrl(P)$. Clearly, $d_1^1d_2^0x = d_1^0d_1^1x \notin Extrl(P)$ and therefore $Extrl(P) \subseteq R$. It follows now from \ref{FBG} that $\iota$ induces an isomorphism of fundamental bipartite graphs and a homotopy equivalence of d-path spaces for each pair of extremal elements. 
\findem

\begin{figure}
\begin{center}
\subfigure[$P$]
{ 
\includegraphics[height=2.5cm]{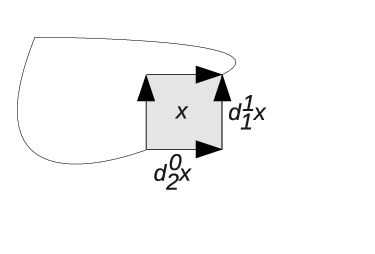}
\label{thm3P}
}
\hspace{2cm}
\subfigure[$Q$]
{
\includegraphics[height=2.5cm]{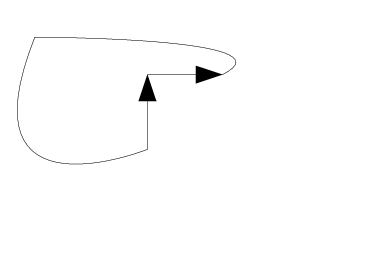}
\label{thm3Q}
}
\caption{Precubical sets $P$ and $Q$ in theorem \ref{elimination1} with $b=0$} \label{thm3}
\end{center}
\end{figure}

The statement of the following theorem is illustrated in figure \ref{thm3}.

\begin{theor} \label{elimination1}
Let $P$ be a precubical set, $b \in \{0,1\}$, and $x \in P_2$ be a regular element such that no element of $P_2\setminus \{x\}$ has $d_1^{1-b}x$  or $d_2^bx$ in its boundary and the only elements of $P_1$ having $d_1^{1-b}d_2^bx = d_1^bd_1^{1-b}x$ in their boundary are $d_1^{1-b}x$ and  $d_2^bx$. Then a precubical subset $Q$ of $P$ such that the inclusion $\iota\colon |Q| \hookrightarrow |P|$ is a dihomotopy equivalence relative to $|Q|$ is given by $Q_0 = P_0\setminus \{d_1^{1-b}d_2^bx\}$, $Q_1 = P_1\setminus \{d_1^{1-b}x, d_2^bx\}$, $Q_2 = P_2\setminus \{x\}$, and $Q_r = P_r$ $(r > 2)$. Moreover, $Extrl(P) = Extrl(Q)$ and $\iota$ induces an isomorphism of fundamental bipartite graphs and a homotopy equivalence of d-path spaces for each pair of extremal elements.
\end{theor}

\dem
The precubical set $P$, the number $b$, and the element $x$ satisfy  the conditions of \ref{elimination2} with $a=1$. It follows that a precubical subset $M$ of $P$ such that $Q$ is a precubical subset of $M$ and the inclusion $|M| \hookrightarrow |P|$ is a dihomotopy equivalence relative to $|Q|$ is given by $M_0 = P_0$, $M_1 = P_1\setminus \{d_{2}^bx\}$, $M_2 = P_2\setminus \{x\}$, and $M_r = P_r$ $(r >2)$. The precubical set $M$, the number $1-b$, and the element $d_1^{1-b}x$ satisfy the conditions of \ref{elimination3}. We have $Y = \{y \in M_1\, |\, d_1^ {1-b}y = d_1^ bd_1^ {1-b}x\} = \emptyset$. It follows that the inclusion $|Q| \hookrightarrow |M|$ is a dihomotopy equivalence relative to $|Q|$ and hence that the inclusion $\iota\colon |Q| \hookrightarrow |P|$ is a dihomotopy equivalence relative to $|Q|$. Note that $d_1^{1-b}d_2^bx = d_1^bd_1^{1-b}x \notin Extrl(P)$. Therefore $Extrl(P) \subseteq Q$. This implies that $Extrl(P) \subseteq Extrl(Q)$. Let $v \in Q_0$ be a minimal element of $Q$. Then for all $z \in Q_1$, $d_1^1z \not= v$. If $b=0$, then $d_2^1x \in Q_1$ and therefore $d_1^1d_1^{1-b}x = d_1^1d_2^1x \not= v$. If $b= 1$, then $d_1^1d_1^{1-b}x = d_1^{1-b}d_2^bx \notin Q_0$ and therefore $d_1^1d_1^{1-b}x \not= v$.  If $b=0$, then $d_1^1d_2^bx = d_1^{1-b}d_2^bx \notin Q_0$ and therefore $d_1^1d_2^bx \not= v$. If $b=1$, then $d_1^1x \in Q_1$ and therefore $d_1^1d_2^bx = d_1^1d_1^1x \not=v$. It follows that $v$ is a minimal element of $P$. A similar argument shows that any maximal element of $Q$ is also a maximal element of $P$. It follows that $Extrl(P) = Extrl(Q)$ and hence, by \ref{FBG}, that $\iota$ induces an isomorphism of fundamental bipartite graphs and a homotopy equivalence of d-path spaces for each pair of extremal elements.  
\findem

\begin{rem}\rm
Note that the equivalence between $|P|$ and $|Q|$ in \ref{elimination1} can be seen as a composition of a T- and a S-homotopy equivalence in the sense of \cite{GaucherGoubault}.
\end{rem}

\section{Examples}

In this section, we use our reduction techniques to compute small models of three well-known example precubical sets. In each case, we know \emph{a priori} that the geometric realizations of the model and the given precubical set are dihomotopy equivalent relative to the extremal elements and have isomorphic fundamental bipartite graphs and homotopy equivalent d-path spaces between extremal elements. In order to construct the small models, we use basically the following straightforward procedure: We run through the 2-dimensional cubes as long as it is possible to eliminate one and after that we proceed similarly with the 1-dimensional cubes.

\begin{ex}\rm \label{exholes}
Consider the $2$-dimensional precubical set depicted in figure \ref{holes} below. The grey squares represent the elements of degree $2$, the arrows represent the elements of degree $1$, and the end points of the arrows represent the elements of degree $0$. The arrow corresponding to an edge $x$ points from $d_1^0x$ to $d_1^1x$. The left-hand edge of a square $x$ is $d_1^0x$, the right-hand edge is $d_1^1x$, the lower edge is $d_2^0x$, and the upper edge is $d_2^1x$.
We use the following sequence of $2$-dimensional reductions to deform this precubical set into the $1$-dimensional precubical subset depicted in figure \ref{holes1}: We proceed linewise from the top left square to the bottom right square using Theorem \ref{elimination1} with $b=1$ to eliminate all squares except for the four squares on the right of the holes where we use Theorem \ref{elimination2} with $a=2$ and $b=0$ and the four squares below the holes where we use Theorem \ref{elimination2} with $a=1$ and $b=1$. A sequence of $1$-dimensional reductions using Theorem \ref{elimination3} with $b=0$ permits us to simplify the model further to the $1$-dimensional precubical set with four vertices in figure \ref{holes2}. From this we obtain the final model in figure \ref{holes3} using Theorem \ref{elimination3} twice with $b=1$. 
\end{ex}

\begin{figure}[H] 
\begin{center}
\subfigure[]
{ 
\includegraphics[height=2.5cm]{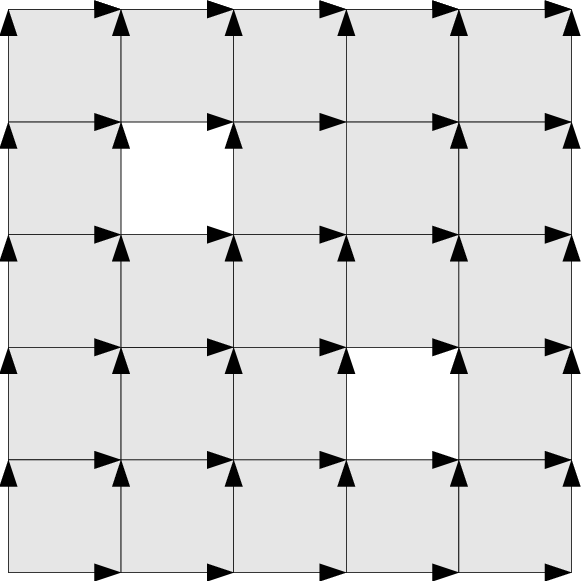}
\label{holes}
}
\subfigure[]
{
\includegraphics[height=2.5cm]{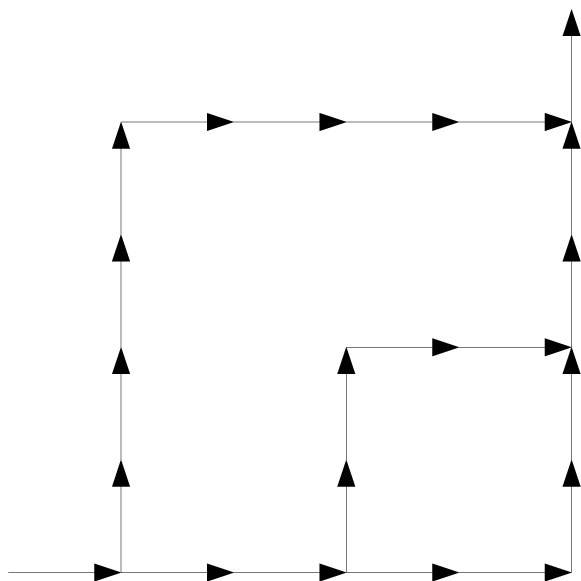}
\label{holes1}
}
\subfigure[]
{
\includegraphics[height=2.5cm]{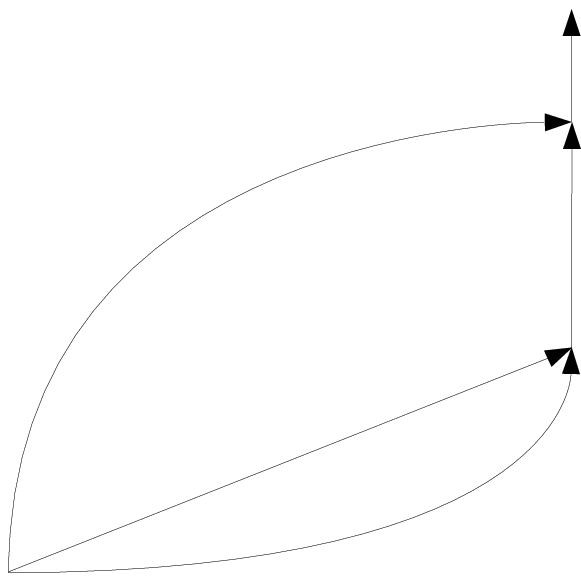}
\label{holes2}
}
\subfigure[]
{
\includegraphics[height=2.5cm]{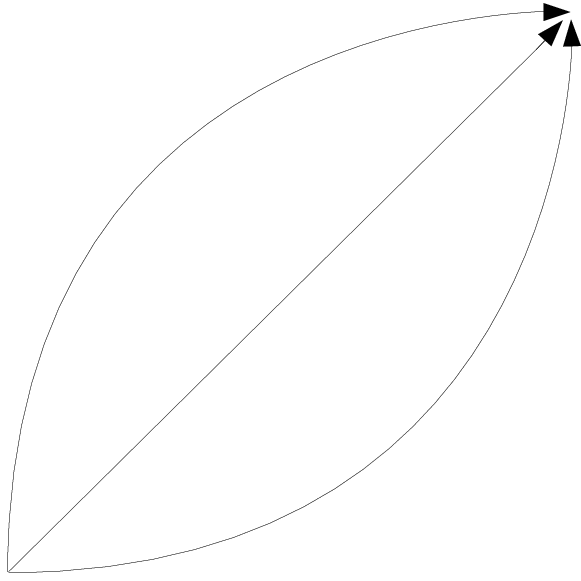}
\label{holes3}
}
\caption{Reduction of a square with two unordered holes}
\end{center}
\end{figure}

\begin{figure}[H]
\begin{center}
\subfigure[]
{ 
\includegraphics[height=2.5cm]{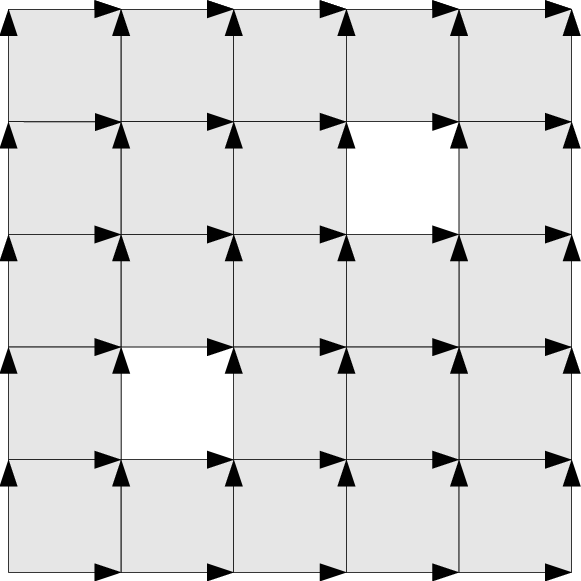}
\label{ordholes}
}
\subfigure[]
{
\includegraphics[height=2.5cm]{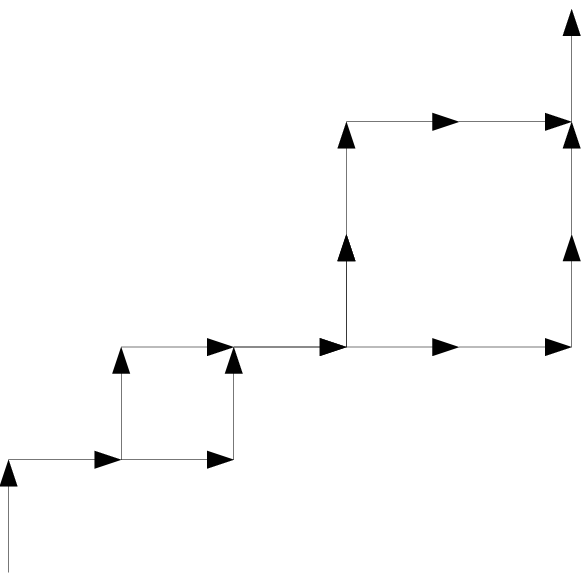}
\label{ordholes1}
}
\subfigure[]
{
\includegraphics[height=2.5cm]{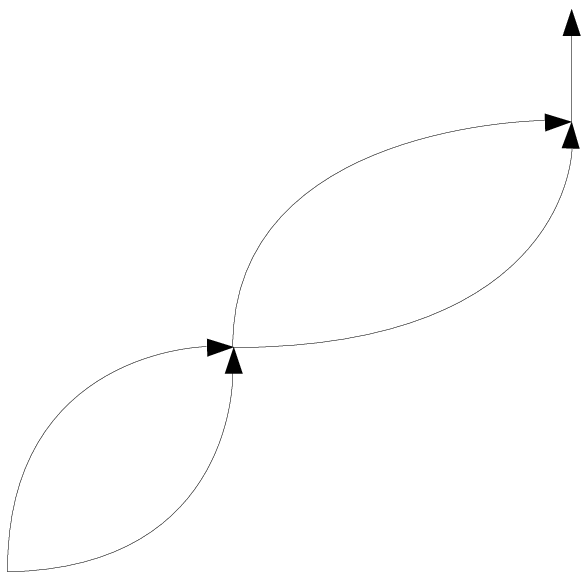}
\label{ordholes2}
}
\subfigure[]
{
\includegraphics[height=2.5cm]{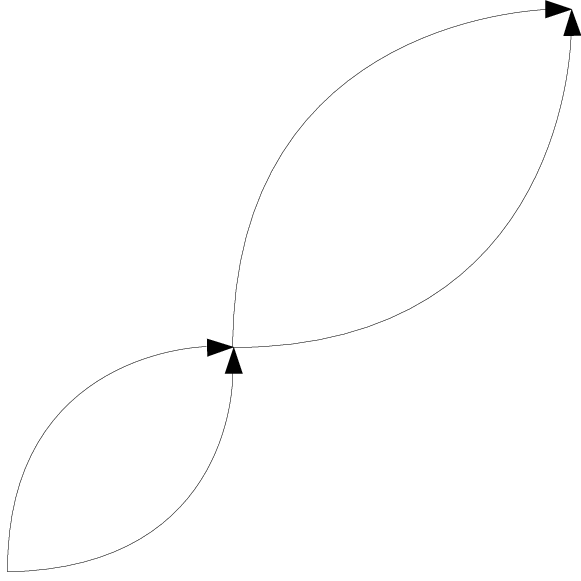}
\label{ordholes3}
}
\caption{Reduction of a square with two ordered holes}
\end{center}
\end{figure}

\begin{ex}\rm \label{exordholes}
Consider the precubical set in figure \ref{ordholes}. We use the following sequence of $2$-dimensional reductions to deform this precubical set into the $1$-dimensional precubical subset in figure \ref{ordholes1}: We proceed linewise from the top left square to the square on the left of the lower hole using Theorem \ref{elimination1} with $b=1$ to eliminate all squares except for the one on the right of the upper hole where we use Theorem \ref{elimination2} with $a=2$ and $b=0$ and the one below the upper hole where we use Theorem \ref{elimination2} with $a=1$ and $b=1$. We then eliminate the remaining squares using Theorem \ref{elimination1} with $b=0$ proceeding linewise upwards from the bottom right square to the square to the right of the lower hole. We simplify the model further to the precubical set in figure \ref{ordholes2} by means of a sequence of $1$-dimensional reductions using Theorem \ref{elimination3} with $b=0$. Using Theorem \ref{elimination3} with $b=1$ we finally obtain the model in figure \ref{ordholes3}.
\end{ex}

\begin{figure}[H]
\begin{center}
\subfigure[]
{ 
\includegraphics[height=2.5cm]{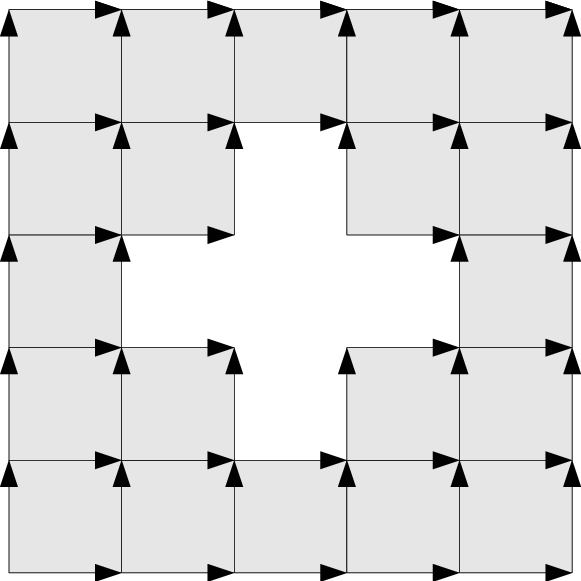}
\label{swiss}
}
\subfigure[]
{
\includegraphics[height=2.5cm]{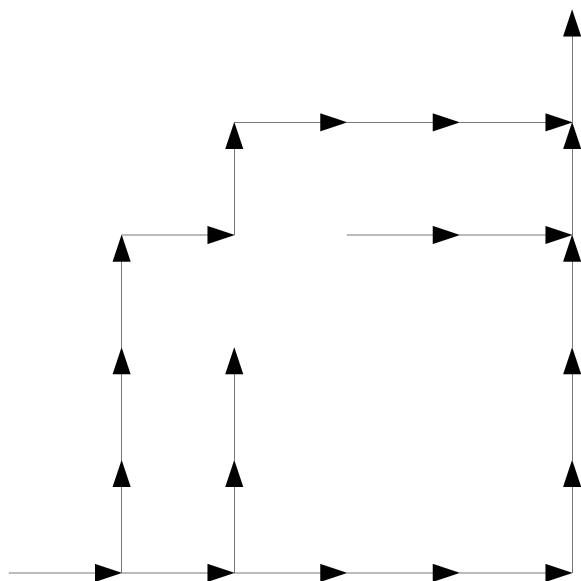}
\label{swiss1}
}
\subfigure[]
{
\includegraphics[height=2.5cm]{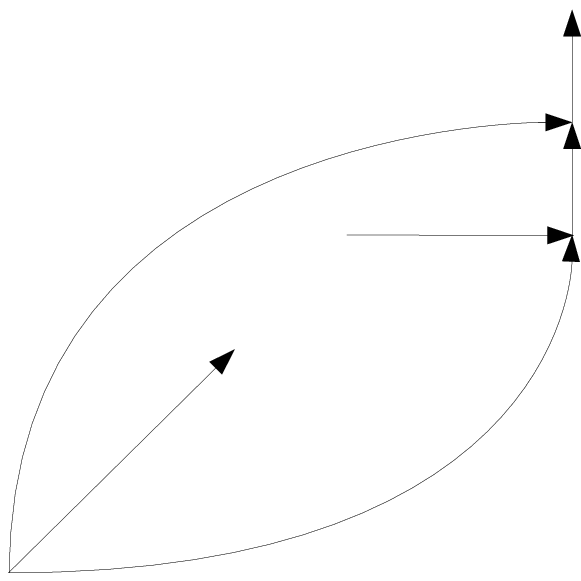}
\label{swiss2}
}
\subfigure[]
{
\includegraphics[height=2.5cm]{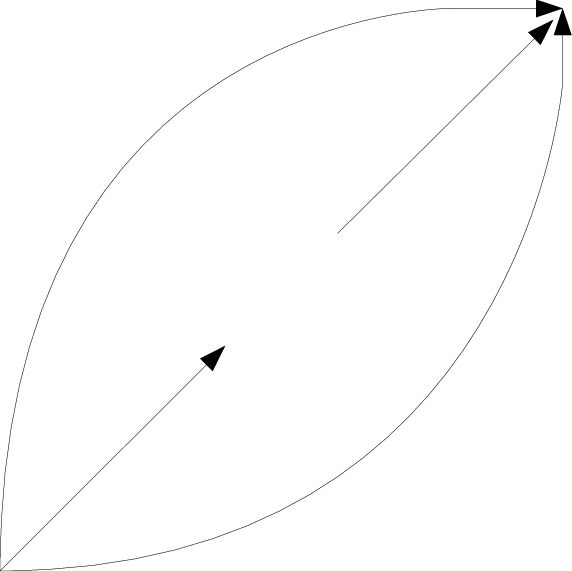}
\label{swissred4}
}
\caption{Reduction of the Swiss flag}
\end{center}
\end{figure}

\begin{ex}\rm \label{exswiss}
We deform the \emph{Swiss flag} in figure \ref{swiss} into the $1$-dimensional precubical subset in figure \ref{swiss1} successively as follows: We proceed linewise from the top left square to the bottom right square using Theorem \ref{elimination1} with $b=1$ to eliminate all squares except for the three squares on the right of the upper and middle holes where we use Theorem \ref{elimination2} with $a=2$ and $b=0$ and the three squares below the left and the middle holes where we use Theorem \ref{elimination2} with $a=1$ and $b=1$. Using Theorem \ref{elimination3} several times with $b=0$ we obtain the $1$-dimensional precubical set in figure \ref{swiss2}. We finally obtain the model in figure \ref{swissred4} using Theorem \ref{elimination3} twice with $b=1$. The reader might be interested to compare this model to the one obtained in \cite{GaucherGoubault} using a sequence of S- and T-homotopy equivalences.
\end{ex}

\section{Final remarks}

\subsection{Further collapsing operations} 

\begin{figure}[H]
\begin{center}
\includegraphics[height=2.5cm]{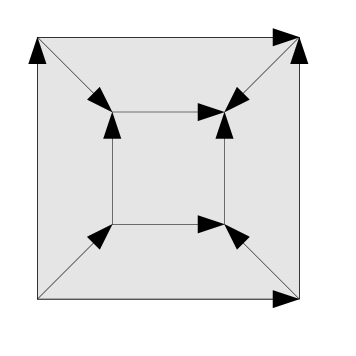}
\caption{Surface of the 3-cube without the bottom face} \label{last}
\end{center}
\end{figure}

\begin{figure}
\begin{center}
\subfigure[$P$]
{ 
\includegraphics[height=2.5cm]{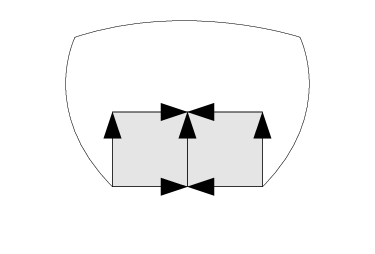}
\label{thm4P}
}
\hspace{2cm}
\subfigure[$Q$]
{
\includegraphics[height=2.5cm]{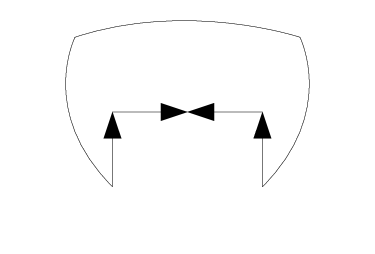}
\label{thm4Q}
}
\caption{$|P|$ and $|Q|$ are dihomotopy equivalent relative to $|Q|$} \label{thm4}
\end{center}
\end{figure}

In this paper, we have established some results which give local combinatorial conditions for the collapsibility of cubes in a 2-dimensional precubical set. The list of our collapsing operations is not exhaustive, and it is possible and for some purposes necessary to establish further collapsibility criteria. Consider, for example, the precubical set depicted in figure \ref{last}, which represents the surface of the 3-cube without the bottom face.
\noindent The geometric realization of this precubical set is dihomotopy equivalent relative to the extremal elements to the directed interval $\vec I$. The results of this paper do not permit us to establish this equivalence. Using the homotopy $H$ defined in the proof of theorem \ref{elimination2} and depicted in figure \ref{H} simultaneously on two adjacent squares, it is, however, easy to establish that the geometric realizations of two precubical sets $P$ and $Q$ that look like the ones in figure \ref{thm4} are dihomotopy equivalent relative to $|Q|$. This result together with the results of this paper permits us to reduce the precubical set of figure \ref{last} to an edge.

\subsection{Higher dimensions} 
\begin{figure}[H]
\begin{center}
\subfigure[Precubical 3-cube ${\mathbb I}^3$]
{ 
\includegraphics[height=3cm]{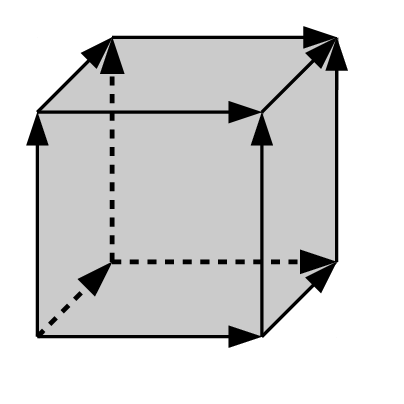}
\label{fullcube}
}
\hspace{2cm}
\subfigure[Precubical subset $Q$]
{
\includegraphics[height=3cm]{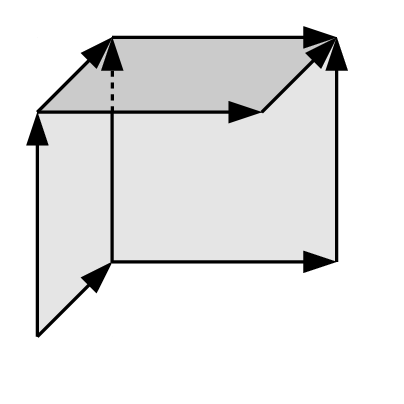}
\label{cubeQ}
}
\caption{$|{\mathbb I}^3|$ and $|Q|$ are not dihomotopy equivalent relative to $|Q|$} \label{cube}
\end{center}
\end{figure}

It is natural to ask whether the  results of this paper can easily be extended to higher dimensions. Unfortunately, the situation in higher dimensions is more complicated. Consider, for example, the precubical 3-cube ${\mathbb I}^3$ and the precubical subset $Q$ depicted in figure \ref{cube}. It is not difficult to see that $|{\mathbb I}^3|$ and $|Q|$ are \emph{not} dihomotopy equivalent relative to $|Q|$. Therefore the three-dimensional version of theorem \ref{elimination1} is not true. Higher dimensional collapsibility results will have either weaker conclusions or stronger and more complex conditions than the results of the present paper. What can be done in higher dimensions is currently being worked out and will be discussed in a forthcoming paper.

\subsection{Extremal models} In \cite{BubenikExtremal}, P. Bubenik introduces \emph{extremal models} of d-spaces and calculates such extremal models for the geometric realizations of the precubical set of the introduction and the precubical sets of Examples \ref{exordholes} and \ref{exswiss}. In all cases, the extremal model is a full subcategory of the fundamental category of the geometric realization of our small model, namely the full subcategory generated by the vertices of the model. It would be interesting to know whether this link between the models constructed using our approach and the extremal models of \cite{BubenikExtremal} can be established in general.


\end{document}